\documentclass[12pt]{article}
\usepackage{amssymb,amsfonts,amsthm}
\usepackage{amstext}

\usepackage{amsmath}

\usepackage{color}
\usepackage[latin1]{inputenc}
\def\openC{{\rm C\kern-.18cm\vrule width.8pt height 7pt depth-.2pt \kern.18cm}}
\def\openN{{{\rm I}\kern-.16em {\rm N}}}
\def\openR{{{\rm I}\kern-.16em {\rm R}}}
\def\openT{{{\rm T}\kern-.42em {\rm T}}}
\def\openZ{{{\rm Z}\kern-.28em{\rm Z}}}

\def\eop{\hfill\rule{2.5mm}{2.5mm}}
\def\pf{\par\smallbreak\noindent {\bf Proof.} }

\newtheorem{thm}{Theorem}[section]

\newtheorem{lem}[thm]{Lemma}
\newtheorem{prop}[thm]{Proposition}

\theoremstyle{definition}

\def\eop{\hfill\rule{2.5mm}{2.5mm}}

\begin{document}

\title{\textbf{An extension of a theorem of Schoenberg to products of spheres} \vspace{-4pt}
\author{\sc
J. C. Guella\thanks{All authors partially supported by FAPESP under grants $\#$2012/22161-3 , $\#$2014/00277-5 and $\#$2014/25796-5 respectively.}, V. A. Menegatto and A. P. Peron}}
\date{}
\maketitle \vspace{-30pt}
\bigskip

\begin{center}
\parbox{13 cm}{{\small We present a characterization for the continuous, isotropic and positive definite kernels on a product of spheres along the lines of a classical result of I. J. Schoenberg on positive definiteness on a single sphere.\ We also discuss a few issues regarding the characterization, including topics for future investigation. }}
\end{center}
\vspace*{1cm}
\noindent{\bf Mathematics Subject Classifications (2010):} 43A35, 33C50, 33C55, 42A10, 42A82\\
\\
\noindent{\bf Keywords:} Positive definiteness, isotropy, spheres, spherical harmonics, Gegenbauer polynomials, addition formula.

\thispagestyle{empty}

%
%

\section{Introduction}\label{s1}

We consider the problem of characterizing positive definite kernels on a product of spheres.\ The focus will be on continuous and isotropic kernels, keeping the setting originally adopted by I. J. Schoenberg in his influential paper published in 1942 (\cite{schoen}).

As usual, let $S^{m}$ denote the unit sphere in the $(m+1)$-dimensional space $\mathbb{R}^{m+1}$ and $S^\infty$ the unit sphere in $\mathbb{R}^\infty$, the usual real $\ell^2$ space.\ Throughout the paper, we will be dealing with real, continuous and isotropic kernels on the product $S^{m}\times S^{M}$, $m,M=1,2,\ldots, \infty$.\ When speaking of continuity, we will assume each sphere is endowed with its usual geodesic distance.\ The {\em isotropy} (zonality) of a kernel $K$ on $S^{m} \times S^{M}$ refers to the fact that
$$
K((x,z),(y,w))=f(x\cdot y,z \cdot w), \quad x,y \in S^{m},\quad z,w \in S^{M},
$$
for some real function $f$ on $[-1,1]^2$, where $\cdot$ stands for the inner product of both $\mathbb{R}^{m+1}$ and $\mathbb{R}^{M+1}$.\ In particular, the concept introduced above demands the usual notion of isotropy on each sphere involved.\ In many places in the paper, we will refer to $f$ as the {\em isotropic part} of $K$.

Recall that if $X$ is a nonempty set, a kernel $K$ is {\em positive definite} on $X$ if
$$
\sum_{\mu,\nu=1}^n c_\mu c_\nu K(x_\mu, x_\nu) \geq 0,
$$
for $n\geq 1$, distinct points on $X$, and reals scalars $c_1, c_2, \ldots, c_n$.\ In other words, for any $n\geq 1$ and any distinct points $x_1, x_2, \ldots, x_n$ on $X$, the $n\times n$ matrix with entries $K(x_\mu,x_\nu)$ is nonnegative definite.\ In this paper, we will present a characterization for the positive definiteness of a continuous and isotropic kernel on $X=S^m \times S^M$ based upon Fourier expansions.

Isotropy and positive definiteness for kernels on a single sphere were first considered by I. J. Schoenberg in \cite{schoen}.\ He showed that a continuous and isotropic kernel $K$ on $S^{m}$ is positive definite if and only if $K(x,y)=g(x\cdot y)$, $x,y \in S^m$, in which the isotropic part $g$ of $K$ has a series representation in the form
$$
g(t)=\sum_{k=0}^{\infty} a_k^m P_k^{m}(t),\quad t \in [-1,1],
$$
in which $a_k^m \geq 0$, $k \in \mathbb{Z}_+$ and $\sum_{k=0}^{\infty}a_k P_k^{m}(1) <\infty$.\ The symbol $P_k^{m}$ stands for the usual
Gegenbauer polynomial of degree $k$ associated with the rational $(m-1)/2$, as discussed in \cite{szego}.\ This Schoenberg's outstanding result is far-reaching and has ramifications in distance geometry, statistics, spherical designs, approximation theory, etc.\ In approximation theory, positive definite kernels are used in interpolation of scattered data over the sphere.\ The importance of this problem in many areas of science and engineering is reflected in the literature, where different methods to solve such problem have been proposed.\ Given $n$ distinct data points $x_1, x_2, \ldots, x_n$ on $S^{m}$ and a target
function $h:S^{m} \to \mathbb{R}$,
the interpolation problem itself requires the finding of a continuous function $s:S^{m} \to \mathbb{R}$ of the form
$$
s(x)=\sum_{j=1}^n \lambda_j g(x \cdot x_j),\quad x \in S^{m},\quad \lambda_1, \lambda_2, \ldots, \lambda_n \in \mathbb{R},$$
so that $s(x_i)=h(x_i)$, $i=1,2,\ldots,n$.
If we choose the prescribed function $g$ to be the isotropic part of a convenient positive definite kernel $K$, then the interpolation problem has a unique solution for any $n$ and any $n$ data points.

An outline of the paper is as follows.\ In Section 2, we present several technical results that culminate with a characterization for the continuous, isotropic and positive definite kernels on $S^m \times S^M$, $m,M <\infty$.\ In Section 3, we complete this circle of ideas, by reaching a similar characterization in the cases in which at least one of the spheres involved is the real Hilbert sphere $S^\infty$.\ Finally, Section 4 contains a few relevant remarks along with the description of future lines of investigation on the subject.

\section{Positive definiteness on $S^{m}\times S^{M}$, $m,M < \infty$.}\label{s2}

This section contains all the technical material needed in the proof of the extension of Schoenberg's theorem to a product of spheres $S^m \times S^M$, in the case when both $m$ and $M$ are finite.\ The proof will require a series of well-known results involving Gegenbauer polynomials and also a few facts from the analysis on the sphere.\ We suggest the classical reference \cite{szego} for the first topic and \cite{atkinson,dai,groemer,muller} for the other.

The orthogonality relation for Gegenbauer polynomials reads as follows (\cite[p.10]{dai}):
$$
\int_{-1}^1 P_n^{m}(t) P_k^{m}(t) (1-t^2)^{(m-2)/2}dt=\frac{\tau_{m+1}}{\tau_{m}}\frac{m-1}{2n+m-1}P_n^{m}(1)\delta_{n,k},
$$
in which $\tau_{m+1}$ is the surface area of $S^{m}$, that is,
$$
\tau_{m+1}:=\frac{2\pi^{(m+1)/2}}{\Gamma((m+1)/2)}.
$$

Since Schoenberg's characterization for positive definiteness on $S^{m}$ is based upon Fourier expansions with respect to the orthogonal family $\{P_n^{m}:n=0,1,\ldots\}$, it is quite natural to expect that a similar characterization for positive definiteness on $S^{m}\times S^{M}$ will require expansions with respect to the tensor family
$$\{(t,s) \in [-1,1]^2 \to P_k^m(t)P_l^M(s): k,l=0,1,\ldots\}.$$

The first important fact to be noticed about the functions in the family above is this.

\begin{lem} \label{gegpd}
If $k,l \in \mathbb{Z}_+$, then $(t,s) \in [-1,1]^2 \to P_k^m(t)P_l^M(s)$ is the isotropic part of a positive definite kernel on $S^{m}\times S^{M}$.
\end{lem}
\pf This follows from the definition of positive definiteness, Schoenberg's original characterization for positive definite kernels and the Schur product theorem (\cite[p.458]{horn}).\ The later asserts that the entry-wise product of two nonnegative definite matrices of same order is a nonnegative definite matrix itself.\eop

The tensor family is orthogonal on $[-1,1]^2$ with respect to the weight function
$$
w_{m,M}(t,s)=(1-t^2)^{(m-2)/2}(1-s^2)^{(M-2)/2},\quad t,s \in [-1,1].
$$
The $(k,l)$-Fourier coefficient of a function $f: [-1,1]^2 \to \mathbb{R}$ from $L^1([-1,1]^2, w_{m,M})$
is
$$
\hat{f}_{k,l}:=\frac{1}{\tau_k^m \tau_l^M}\int_{[-1,1]^2} f(t,s)P_k^{m}(t)P_l^{M}(s)dw_{m,M}(t,s), \quad k,l \in \mathbb{Z}_+,
$$in which
$$
\tau_k^m:=\frac{\tau_{m+1}}{\tau_{m}}\frac{m-1}{2k+m-1}P_k^{m}(1),\quad k\in \mathbb{Z}_+.
$$
The next lemma describes an alternative way for computing these Fourier coefficients.\ The symbol $\sigma_m$ will denote the surface measure on $S^m$.

\begin{lem} \label{pdmultiple} If $f$ belongs to $L^1([-1,1]^2,w_{m,M})$, then the Fourier coefficient $\hat{f}_{k,l}$ is a positive constant multiple of
$$
\int_{S^{m}\times S^{M}} \left[\int_{S^{m}\times S^{M}} f(x\cdot y,z \cdot w)P_k^{m}(x\cdot y) \times P_l^{M}(z\cdot w)d\sigma_{m}(y)d\sigma_{M}(w)\right]d\sigma_{m}(x)d\sigma_{M}(z).
$$
\end{lem}
\pf If $m,M \geq 2$, it suffices to employ the Funk-Hecke formula (\cite[p.11]{dai}) in the expression defining the Fourier coefficient.\
The Funk-Hecke formula states that
$$
\int_{S^{M}}g(z\cdot w)P_l^{M}(z\cdot w)d\sigma_{M}(w)=\tau_{m-1}\int_{-1}^1 g(s)P_l^{M}(s)(1-s^2)^{(M-2)/2}ds,\quad z \in S^{M},
$$
whenever $l\in \mathbb{Z}_+$ and $g \in L^2([-1,1],w_M)$.\ Using the formula with
$$
g(s)=\int_{-1}^1 f(t,s)P_k^{m}(t)(1-t^2)^{(m-2)/2}dt, \quad s \in [-1,1],
$$
it is promptly seen that $\hat{f}_{k,l}$ is a positive multiple of
$$
\int_{S^{M}} \int_{-1}^1 f(t,z \cdot w)P_k^{m}(t)(1-t^2)^{(m-2)/2}dt \ P_l^{M}(z\cdot w) d\sigma_{M}(w).
$$
Applying a similar argument in the internal integral reveals that $\hat{f}_{k,l}$ is a positive multiple of
$$
\int_{S^{M}} \int_{S^{m}} f(x\cdot y,z \cdot w)P_k^{m}(x\cdot y) P_l^{M}(z\cdot w)d\sigma_{m}(y)d\sigma_{M}(w).
$$
Integration with respect to the remaining variables concludes the proof.\ In the cases in which either $m=1$ or $M=1$, the arguments demand the replacement of the Funk-Hecke formula with direct computation.\eop

\begin{lem} \label{pdtimes} If $f$ is the continuous and isotropic part of a positive definite kernel on $S^{m}\times S^{M}$, then
$$
\int_{S^{m}\times S^{M}}\left[\int_{S^{m}\times S^{M}}f(x\cdot y,z\cdot w)d\sigma_{m}(x)d\sigma_{M}(z)\right]d\sigma_{m}(y)d\sigma_{M}(w) \geq 0.
$$
\end{lem}
\pf It suffices to write the double integral $I$ in the statement of the theorem as a double limit of Riemann sums.\ Indeed, we can select a sequence $\{\mathcal{P}_n: n=0,1,\ldots\}$ of partitions of $S^{m}\times S^{M}$
in such a way that $\mathcal{P}_n=\{Q_1^n,Q_2^n,\ldots, Q_{\alpha(n)}^n\}$, the sequence $\{\alpha(n)\}$ increases to $\infty$ and the sequences of diameters $\{\mbox{diam}(Q_j^n)\}$ satisfy
$\lim_{n\to \infty}\mbox{diam}(Q_j^n)=0$.\ Picking points $(x_j^n,z_j^n) \in Q_j^n$, we can write
$$
I=\lim_{N \to \infty} \sum_{J=1}^{\alpha(N)}\left[ \int_{S^{m}\times S^{M}}f(x \cdot x_J^N,z \cdot z_J^N)d\sigma_{m}(x)d\sigma_{M}(z)\right]\mbox{vol}(Q_J^N).
$$
Repeating the procedure with the resulting integral leads to
\begin{eqnarray*}
I & = & \lim_{N \to \infty} \sum_{J=1}^{\alpha(N)}\left[\lim_{n \to \infty}\sum_{j=1}^{\alpha(n)}f(x_j^n \cdot x_J^N,z_j^n \cdot z_J^N) \mbox{vol}(Q_j^n)\right] \mbox{vol}(Q_J^N)\\
  & = & \lim_{N \to \infty} \lim_{n \to \infty}\sum_{J=1}^{\alpha(N)}\sum_{j=1}^{\alpha(n)}\mbox{vol}(Q_j^n)\mbox{vol}(Q_J^N)f(x_j^n \cdot x_J^N,z_j^n \cdot z_J^N).
  \end{eqnarray*}
Since the double limit above exists, it follows that
$$
I=\lim_{n \to \infty}\sum_{j,J=1}^{\alpha(n)}\mbox{vol}(Q_j^n)\mbox{vol}(Q_J^n)f(x_j^n \cdot x_J^n,z_j^n \cdot z_J^n).
$$
If $f$ is the isotropic part of a positive definite kernel on $S^{m}\times S^{M}$, then each double sum in the last expression above is clearly nonnegative.\ In particular, the limit itself is nonnegative as well.\eop

We now combine the three lemmas above in order to obtain the following result.

\begin{lem} \label{coeffpos} If $f$ is the continuous and isotropic part of a positive definite kernel on $S^{m}\times S^{M}$ then
$$
\hat{f}_{k,l} \geq 0, \quad k,l \in \mathbb{Z}_+.
$$
\end{lem}
\pf Let us fix $k$ and $l$.\ Lemma \ref{gegpd} and the Schur product theorem guarantees that the function
$$
(t,s) \in [-1,1]^2 \to f(t,s)P_k^{m}(t)P_l^{M}(s)
$$
is the continuous and isotropic part of a positive definite kernel on $S^{m} \times S^{M}$.\ Taking into account this information and that one provided by Lemma \ref{pdmultiple}, an application of Lemma \ref{pdtimes} leads to the inequality in the statement of the lemma.\eop

Next, we recall one of the several generating formulas for the Gegenbauer polynomials, the Poisson identity (\cite[p.419]{dai}).

\begin{lem} \label{poisson} If $r \in [0,1)$, then
$$
\frac{1-r^2}{(1-2tr+r^2)^{(m+1)/2}}=\sum_{k=0}^{\infty}\frac{2k+m-1}{m-1}P_k^{m}(t)r^k, \quad t\in [-1,1].
$$
If $r_0 \in [0,1)$ is fixed, then the convergence of the series is absolute and uniform for $(r,t) \in [0,r_0]\times [-1,1]$.
\end{lem}

We are about ready to prove the following auxiliary result.

\begin{lem} \label{spoisson} Let $f$ be the continuous and isotropic part of a kernel on $S^{m}\times S^{M}$.\ If $r,\rho \in [0,1)$, then the double series
$$
\sum_{k,l=0}^{\infty}\hat{f}_{k,l}P_{k}^{m}(1)P_{l}^{M}(1)r^k\rho^l
$$
converges.\ As a matter of fact, there exists a positive constant $C$, depending upon $f$ only, so that
$$
\left|\sum_{k,l=0}^{\infty}\hat{f}_{k,l}P_{k}^{m}(1)P_{l}^{M}(1)r^k\rho^l\right|\leq C, \quad r,\rho \in [0,1).
$$
\end{lem}
\pf Let $a_{k,l}^{r,\rho}$ denote the general term of the series in the statement of the lemma.\ It is promptly seen that
$$
a_{k,l}^{r,\rho}=C_1\int_{-1}^1 \int_{-1}^1 f(t,s)\frac{2k+m-1}{m-1}P_k^{m}(t)r^k\frac{2l+M-1}{M-1}P_l^{M}(t)\rho^l dw_{m,M}(t,s),
$$
in which
$$
C_1=\frac{\tau_{m}\tau_{M}}{\tau_{m+1}\tau_{M+1}}.
$$
On the other hand, Lemma \ref{poisson} implies that
$$
\sum_{k=0}^{\infty}\sum_{l=0}^{\infty}a_{k,l}^{r,\rho}=C_1 \int_{-1}^1 \int_{-1}^1f(t,s)\frac{1-r^2}{(1-2rt+r^2)^{(m+1)/2}}\frac{1-\rho^2}{(1-2\rho s+s^2)^{(M+1)/2}}dw_{m,M}(t,s),
$$
whenever $r,\rho \in [0,1)$.\ Thus, since $f$ is continuous and the left hand side of the Poisson identity is positive, the proof of the first half of the lemma reduces itself to proving that the double integral
$$
\int_{-1}^1\int_{-1}^1 \frac{1-r^2}{(1-2rt+r^2)^{(m+1)/2}}\frac{1-\rho^2}{(1-2\rho s+s^2)^{(M+1)/2}}(1-t^2)^{(m-2)/2}(1-s^2)^{(M-2)/2}dtds
$$
is finite.\ But, this follows from the well-known property of the Poisson kernels (\cite[p.47]{muller})
$$
\int_{-1}^1 \frac{1-r^2}{(1-2rt+r^2)^{(m+1)/2}}(1-t^2)^{(m-2)/2}dt=\frac{\tau_{m+1}}{\tau_{m}}.
$$
In particular, $C:=\max\{|f(s,t)|: -1\leq t,s\leq 1\}$ fits into what is needed in the second statement of the lemma.\eop

\begin{lem}\label{pd1} If $f$ is the continuous and isotropic part of a positive definite kernel on $S^{m}\times S^{M}$, then the double series
$$
\sum_{k,l=0}^{\infty}\hat{f}_{k,l}P_{k}^{m}(1)P_{l}^{M}(1)
$$
converges.
\end{lem}
\pf Let $f$ be the continuous and isotropic part of a positive definite kernel on $S^{m}\times S^{M}$.\ Due to Lemma \ref{coeffpos}, we know already that all the Fourier coefficients $\hat{f}_{k,l}$ are nonnegative.\ In particular, the double sequence $\{s_{p,q}\}$ of partial sums of the double series in the statement of the current lemma is monotonically increasing, that is, $s_{p,q} \leq s_{\mu,\nu}$ when $p\leq \mu$ and $q \leq \nu$.\ On the other, the previous lemma produces the inequality
$$
\sum_{k=0}^{p}\sum_{l=0}^{q}\hat{f}_{k,l}P_{k}^{m}(1)P_{l}^{M}(1)r^k\rho^l \leq C, \quad p,q\in\mathbb{Z}_+, \quad r,\rho \in [0,1),
$$
for some $C>0$.\ By taking a double limit when $r,\rho \to 1^{+}$, we deduce that the double sequence of partial sums $\{s_{p,q}\}$ is bounded above.\ A classical result from the theory
of double sequences (\cite[p.373]{limaye}) implies that $\{s_{p,q}\}$ converges, that is, the series in the statement of the lemma converges.\eop

The Weierstrass M-test can be adapted to hold for double series of functions.\ Combining it with Lemma \ref{pd1} leads to the proposition below.

\begin{prop}\label{conve} If $f$ is the continuous and isotropic part of a positive definite kernel on $S^{m}\times S^{M}$, then
$$
\sum_{k,l=0}^{\infty}\hat{f}_{k,l}P_{k}^{m}(t)P_{l}^{M}(s)
$$
converges absolutely and uniformly for $(t,s) \in [-1,1]^2$.
\end{prop}

The main result in this section is as follows.

\begin{thm} \label{mainPD} Let $K$ be a continuous and isotropic kernel on $S^{m}\times S^{M}$.\ It is positive definite on $S^{m}\times S^{M}$ if and only if its isotropic part $f$ has a representation in the form
$$
f(t,s)=\sum_{k,l=0}^{\infty}\hat{f}_{k,l}P_{k}^{m}(t)P_{l}^{M}(s), \quad t,s \in [-1,1],
$$
in which $\hat{f}_{k,l} \geq 0$, $k,l \in \mathbb{Z}_+$ and $\sum_{k,l=0}^{\infty}\hat{f}_{k,l}P_{k}^{m}(1)P_{l}^{M}(1)<\infty$.
\end{thm}
\pf If the isotropic part $f$ of $K$ has the representation announced in the theorem, the series appearing there is uniformly and absolutely convergent.\ In particular, Lemma \ref{gegpd} implies that $f$ is
a pointwise double limit of functions which are isotropic parts of positive definite kernels on $S^{m} \times S^{M}$.\ Consequently, $K$ itself is positive definite on $S^{m}\times S^{M}$.\ Conversely, assume $K$ is positive definite on $S^{m}\times S^{M}$ and write $f$ to denote its isotropic part.\ Lemma \ref{pd1} and Proposition \ref{conve} supply a function $g$ so that
$$
g(s,t)=\sum_{k,l=0}^{\infty}\hat{f}_{k,l}P_{k}^{m}(t)P_{l}^{M}(s), \quad t,s \in [-1,1],
$$
with uniform convergence in $[-1,1]^2$.\ In particular, $g$ is continuous in $[-1,1]^2$.\ On the other hand, the same uniform convergence and the
orthogonality relation mentioned at the beginning of the section imply that
$$
\hat{f}_{k,l} -\hat{g}_{k,l}=0, \quad k,l \in \mathbb{Z}_+.
$$
Consequently, $f=g$.\eop

\section{Positive definiteness on $S^\infty \times S^M$}

In this section, we extend Theorem \ref{mainPD} to the cases in which either $m=\infty$ or $M=\infty$.\ Clearly, it suffices to consider the cases $m=\infty$, $M<\infty$ and $m=M=\infty$ only.

Every sphere $S^{m}$ can be isometrically embedded in $S^{\infty}$.\ In particular, a positive definite kernel on $S^\infty \times S^{M}$ is positive definite on $S^{m}\times S^{M}$, for $m=1,2,\ldots$.\ Likewise, if $f$ is the isotropic part of a positive definite kernel on $S^\infty \times S^{M}$, then it is the isotropic part of a positive definite kernel on $S^{m}\times S^{M}$, for $m=1,2, \ldots$.\ In addition, if $f$ is continuous, then for
every $m\geq 1$, we have a representation for $f$ in the form
$$
f(t,s)=\sum_{k,l=0}^{\infty}\hat{f}_{k,l}^{m,M} P_{k}^{m}(t) P_{l}^{M}(s), \quad t,s \in [-1,1],
$$
in which
$$
\hat{f}_{k,l}^{m,M} =\frac{1}{\tau_{k}^{m} \tau_{l}^{M}}  \int_{[-1,1]^{2}} f(x,y) P_{k}^{m}(t)P_{l}^{M}(s) dw_{m,M}(t,s)\geq 0,\quad k,l \in \mathbb{Z}_+.
$$
and $\sum_{k=0}^{\infty}\sum_{l=0}^\infty \hat{f}_{k,l}^{m,M}P_{k}^{m-1}(1) P_{l}^{M-1}(1)<\infty$.\ Below, we will prefer to normalized the above expressions by writing
$$
R_k^{m}=\frac{P_k^{m}}{P_k^{m}(1)}, \quad k\in \mathbb{Z}_+,
$$
and
$$
f(t,s)=\sum_{k,l=0}^{\infty} \check{f}_{k,l}^{m,M} R_{k}^{m}(t) R_{l}^{M}(s), \quad t,s \in [-1,1],
$$
where now
$$
\check{f}_{k,l}^{m,M}=P_k^{m}(1)P_l^{M}(1)\hat{f}_{k,l}^{m,M}, \quad k,l \in \mathbb{Z}_+.
$$
Before we proceed, it is convenient to mention that the Fourier coefficients introduced above are well-defined as long as $f$ belongs to $L^1([-1,1],w_{m,M})$.

\begin{lem}\label{pre} Let $f$ belong to $L^1([-1,1],w_{m,M})$.\ If $k$ and $l$ are fixed nonnegative integers, then the sequence $\{\check{f}_{k,l}^{2m,M}:m=1,2,\ldots\}$ is convergent.
\end{lem}
\pf Using the following recurrence relation for Gegenbauer polynomials (\cite[p. {84}]{szego})
$$ (1-t^{2})P^{m+2}_{k}(t) = \frac{(k+m-1)(k+m)}{(m-1)(2k+m+1)} P_{k}^{m}(t) - \frac{(k+1)(k+2)}{(m-1)(2k+m+1)}P_{k+2}^{m}(t),
$$
it is easy to deduce that
$$
\check{f}_{k,l}^{m+2,M} = \frac{(k+m-1)(k+m)}{m(2k+m-1)} \check{f}_{k,l}^{m,M} - \frac{(k+1)(k+2)}{m(2k+m+3)} \check{f}_{k+2,l}^{m,M},\quad m\geq 1.
$$
Consequently,
\begin{align*}
|\check{f}_{k,l}^{m+2,M} -  \check{f}_{k,l}^{m,M} | & = \left |\frac{k(k-1)}{m(2k+m-1)}\check{f}_{k,l}^{m,M} - \frac{(k+1)(k+2)}{m(2k+m+3)} \check{f}_{k+2,l}^{m,M} \right | \\
& \leq \left[ \frac{k(k-1)}{m(2k+m-1)} + \frac{(k+1)(k+2)}{m(2k+m+3)}\right]f(1,1),\quad m\geq 1.
\end{align*}
As an obvious consequence, $\{\check{f}_{k,l}^{2m,M}\}$ is a Cauchy sequence of real numbers, therefore, convergent. \eop

\begin{lem}\label{prepar}
If $f$ is the continuous and isotropic part of a positive definite kernel on $S^{\infty} \times S^{M}$, then the double series
$$
\sum_{k,l=0}^{\infty} \check{f}_{k,l}^{m,M}\, t^k R_{l}^{M}(s)
$$
converges for $(t,s) \in (-1,1)^2$, uniformly in $m$.
\end{lem}
\pf In order to see that the series converges in $(-1,1)^2$, for every $m$, it suffices to show that $
\sum_{k=0}^{\infty}\sum_{l=0}^\infty \check{f}_{k,l}^{m,M} |t|^k$ converges.\ Recalling Tonelli's theorem for convergence of double series (\cite[p.384]{limaye}),
that will follow as long as
$\sum_{l=0}^\infty \check{f}_{k,l}^{m,M}$ converges for all $k$ and the iterated series $
\sum_{k=0}^{\infty}\left(\sum_{l=0}^\infty \check{f}_{k,l}^{m,M}\right) |t|^k $ converges.\ But both assertions follow from the inequalities
$$
\sum_{l=0}^\infty \check{f}_{k,l}^{m,M} \leq \sum_{\mu,l=0}^\infty  \check{f}_{\mu,l}^{m,M}= f(1,1), \quad k=0,1,\ldots,
$$
and
$$
\sum_{k=0}^{\infty}\left(\sum_{l=0}^\infty \check{f}_{k,l}^{m,M}\right) |t|^k \leq f(1,1)\sum_{k=0}^{\infty}|t|^k =\frac{f(1,1)}{1-|t|}, \quad t\in (-1,1).
$$
As for the uniform convergence in $m$, it suffices to observe that
$$
\sum_{k,l=0}^{\infty}\check{f}_{k,l}^{m,M}\, t^k R_{l}^{M}(s) \leq f(1,1)\sum_{k=0}^{\infty}|t|^k \quad t,s \in (-1,1),
$$
The proof is complete.\eop

The next lemma is a technical result that can be found proved in \cite{schoen}.

\begin{lem}\label{limitscho}
If $t \in (-1,1)$, then the sequence $\{R_k^m(t)\}$ converges to $t^k$ as $m \to \infty$, uniformly in $k$.
\end{lem}

\begin{thm} \label{maininfty} Let $K$ be a continuous and isotropic kernel on $S^{\infty}\times S^{M}$.\ It is positive definite on $S^{\infty}\times S^{M}$ if and only if its isotropic part $f$ has a representation in the form
$$ f(t,s) = \sum_{k,l=0}^{\infty} \check{f}_{k,l}^{M}t^{k}R_l^{M}(s),$$
in which $\check{f}_{k,l}^{M} \geq 0$, $k,l \in \mathbb{Z}_+$ and $\sum_{k,l=0}^{\infty} \check{f}_{k,l}^{M} < \infty$.
\end{thm}
\pf For each $k$ and $l$, the function $(t,s)\in [-1,1]^2 \to t^{k}R_l^{M}(s)$ is the isotropic part of a positive definite kernel on $S^\infty \times S^{M}$.\ Hence, if $f$ has the representation described in the statement of the theorem, then $K$ is a pointwise limit of positive definite kernels.\ In particular, it is positive definite itself.\ Conversely, assume $K$ is positive definite.\ Without loss of generality, we can assume that $K$ is nonzero.\ Hence, we can assume that its isotropic part $f$ satisfies $f(1,1) >0$.\ Since $f$ is the isotropic part of a positive definite kernel on each product $S^{m} \times S^{M}$, then for each pair $(k,l)$, we may consider the sequence of normalized Fourier coefficients $\{\check{f}_{k,l}^{m,M}\}$.\ Lemma \ref{pre} authenticates the definition
$$
\check{f}_{k,l}^{\, M}:= \lim_{m \to \infty}\check{f}_{k,l}^{2m,M}, \quad k,l=0,1,\ldots.
$$
while Lemma \ref{prepar} guarantees that
$$
\lim_{m\to \infty}\sum_{k,l=0}^{\infty} \check{f}_{k,l}^{2m,M}\, t^k R_{l}^{M}(s)=\sum_{k,l=0}^{\infty}\check{f}_{k,l}^{\, M}\, t^k R_{l}^{M}(s), \quad t,s \in (-1,1).
$$
To proceed, we fix $(t,s) \in (-1,1)^2$ and $\epsilon >0$.\ From the previous limit, we can select $m_0$ so that
$$
\left|\sum_{k,l=0}^{\infty}\check{f}_{k,l}^{2m,M}\, t^k R_{l}^{M}(s)- \sum_{k,l=0}^{\infty} \check{f}_{k,l}^{\, M}\, t^k R_{l}^{M}(s)\right| < \frac{\epsilon}{2}, \quad m\geq m_0.$$
By Lemma \ref{limitscho}, we can select $m_1$ so that
$$|R_k^{2m}(t)-t^k| <\frac{\epsilon}{2f(1,1)}, \quad k=0,1,\ldots, \quad m\geq m_1.$$
It is now clear that
\begin{eqnarray*}
\left|\sum_{k,l=0}^{\infty} \check{f}_{k,l}^{2m,M}\, R_k^{2m}(t)R_{l}^{M}(s)-\sum_{k,l=0}^{\infty} \check{f}_{k,l}^{2m,M}\, t^k R_{l}^{M}(s)\right| & < &  \frac{\epsilon}{2f(1,1)}\sum_{k,l=0}^{\infty} \check{f}_{k,l}^{2m,M} \\
& \leq & \frac{\epsilon}{2}, \quad m\geq m_1.
\end{eqnarray*}
Thus, with the help of an arbitrarily large $m$, we can use both implications above to deduce that
$$
0\leq \left| f(t,s)- \sum_{k,l=0}^{\infty} \check{f}_{k,l}^{\, M}\, t^k R_{l}^{M}(s)\right| < \epsilon.
$$
Hence,
$$
f(t,s)= \sum_{k,l=0}^{\infty} \check{f}_{k,l}^{\, M}\, t^k R_{l}^{M}(s), \quad t,s \in (-1,1).
$$
The coefficients in the representation above are obviously nonnegative.\ If $\sum_{k,l=0}^{\infty} \check{f}_{k,l}^{\, M}$ were not convergent, we could select
a positive integer $N$ so that
$$
\sum_{k,l=0}^{N} \check{f}_{k,l}^{M} \geq 2f(1,1).
$$
Picking a $\tau \in (0,1)$ so that $\tau^{N} > 1/2$, we would reach
$$ f(\tau,1)=\sum_{k,l=0}^{\infty} \check{f}_{k,l}^{M}\tau^{k} \geq  \sum_{k,l=0}^{N} \check{f}_{k,l}^{M}\tau^{k} > f(1,1),$$
a contradiction with the positive definiteness of $f$.\ Having guaranteed the uniform convergence of the series in the representation for $f$ above and invoking the continuity of $f$ in $[-1,1]^2$, we now can let $t,s\to 1^{-}$ and $t,s \to -1^{+}$
in the representation formula in order to conclude that it also holds in $[-1,1]^2$.\eop

A standard adaptation of the arguments used in the proof of the previous theorem is all that is needed in order to deduce the following complement.

\begin{thm}  Let $K$ be a continuous and isotropic kernel on $S^{\infty}\times S^{\infty}$.\ It is positive definite on $S^{\infty}\times S^{\infty}$ if and only if its isotropic part $f$ has a representation in the form
$$ f(t,s) = \sum_{k,l=0}^{\infty}  f_{k,l}t^{k}s^l,$$
in which $f_{k,l} \geq 0$, $k,l \in \mathbb{Z}_+$ and $\sum_{k,l=0}^{\infty} f_{k,l} < \infty$.
\end{thm}

\section{Final remarks}

In view of the characterization for the continuous, isotropic and positive definite kernels on a product of the form $S^m \times S^M$ obtained in the previous section, one may ask what are the other relevant questions regarding that class of kernels.\ We will mention a few of them in this final section of the paper along with some additional elementary results.

Let us begin with the strictly positive definite kernels.\ A continuous, isotropic and positive definite kernel $K$ on $S^m \times S^M$ is {\em strictly positive definite of order $n$} on $S^m \times S^M$ if
its isotropic part $f$ satisfies
$$
\sum_{\mu=1}^n\sum_{\nu=1}^n c_\mu c_\nu f(x_\mu \cdot x_\nu, w_\mu \cdot w_\nu) >0,
$$
whenever the $n$ points $(x_1, w_1), (x_2, w_2), \ldots, (x_n,w_n)$ of $S^m \times S^M$ are distinct and the scalars $c_\mu$ are not all zero.\ So, for a fixed $n$, an interesting question would be to characterize, via the main theorems proved here, the continuous, isotropic and strictly positive definite kernels of order $n$ on $S^m \times S^M$.\ Going one step further, to characterize the continuous, isotropic and strictly positive definite kernels of all orders on $S^m \times S^M$.\ Even in the case of a single sphere, similar characterizations are not available for all fixed $n$ (see \cite{mene0}).

Concerning the problems mentioned above, the intermediate problem to be described below could provide clues to a complete solution.\ A continuous, isotropic and positive definite kernel $K$ on $S^m \times S^M$ is {\em $DC$-strictly positive definite of order $n$} on $S^m \times S^M$ if
its isotropic part $f$ satisfies
$$
\sum_{\mu=1}^n\sum_{\nu=1}^n c_\mu c_\nu f(x_\mu \cdot x_\nu, w_\mu \cdot w_\nu) >0,
$$
whenever the $n$ points $x_1$, $x_2, \ldots, x_n$ of $S^m$ are distinct, the $n$ points $w_1, w_2, \ldots, w_n$ of $S^M$ are distinct and the scalars $c_\mu$ are not all zero.\ Obviously, a strictly positive definite kernel of order $n$ on $S^m \times S^M$ is $DC$-strictly positive definite of order $n$ on $S^m \times S^M$, but not conversely (unless $n=1$).\ Thus, to characterize the continuous, isotropic and positive definite kernels on $S^m \times S^M$ which are $DC$-strictly positive definite of order $n$ on $S^m \times S^M$ would be an interesting problem as well.

A third problem we would like to mention is the description of consistent methods to construct continuous, isotropic and (strictly) positive definite kernels on $S^m \times S^M$.\ In particular, methods based on the description via known classes of continuous, isotropic and (strictly) positive definite kernels on a single sphere.

A very elementary one is this.

\begin{prop} \label{prodi} If $f$ is the continuous and isotropic part of a positive definite kernel on $S^{m}$ and $g$ is the continuous and isotropic part of a positive definite kernel on $S^{M}$, then the function $h$ given by the formula
$$
h(t,s)=f(t)g(s), \quad t,s \in [-1,1],
$$
is the isotropic part of a positive definite kernel on $S^{m}\times S^{M}$.\ Further, if $f$ is the isotropic part of a strictly positive definite kernel of order $n$ on $S^{m}$ (respect., $g$ is the isotropic part of a strictly positive definite kernel of order $n$ on $S^{M}$) and $g(1)>0$ (respect., $f(1)>0$), then $h$ is the isotropic part of a strictly positive definite kernel of order $n$ on $S^{m} \times S^{M}$.
\end{prop}
\pf The first assertion of the theorem is a consequence of the Schur product theorem.\ As for the second one, it follows from Oppenheim's inequality (\cite[p.480]{horn}).\eop

If the intention is a more concrete example, one may employ completely monotonic functions in two variables.\ A continuous function $g: [0,\infty)^2 \to \mathbb{R}$ is \emph{completely monotonic} on $(0,\infty)^2$ if it is $C^\infty$ in $(0,\infty)^2$ and
$$
(-1)^{n_1+n_2} \frac{\partial^{n_1+n_2} g}{\partial u^{n_1}\partial v^{n_2}}(u,v) \geq 0, \quad u,v>0, \quad n_1,n_2\in\mathbb{Z}_+.
$$
It is known that function $g$ as above can be represented in the form
\begin{eqnarray} \label{eq-function CM}
g(u,v) = \int_{[0,\infty)^2}e^{-tu-sv}d\rho(t,s),\quad u,v >0,
\end{eqnarray}
in which $\rho$ is a $\sigma$-additive and nonnegative measure on $[0,\infty)^2$ satisfying $0<\rho((0,\infty)^2)\leq\rho([0,\infty)^2)\leq\infty$ (\cite[p. 87]{bochner}).

A positive scalar multiple of a completely monotonic function on $(0,\infty)^2$ is itself completely monotonic on $(0,\infty)^2$.\ Likewise, the sum and product of two completely monotonic function
on $(0,\infty)^2$ are completely monotonic on $(0,\infty)^2$.\ If $g,h: [0,\infty) \to \mathbb{R}$ are usual completely monotonic functions on $(0,\infty)$, then $F(u,v)=g(u)h(v)$ is completely monotonic
on $(0,\infty)^2$.\ In particular, $(u,v)\in [0,\infty)^2\to \exp(-u)\exp(-v)$
 and $(u,v)\in [0,\infty)^2\to 1/(1+u)^\alpha(1+v)^\beta$, $\alpha,\beta \geq 0$, are completely monotonic on $(0,\infty)^2$.\ Additional examples can be found in \cite{samko}.

For actual examples of positive definite kernels on $S^m \times S^M$, the following result is quite useful.

\begin{prop}
If $g$ is completely monotonic on $(0,\infty)^2$, then
$$
f(t,s):= g(\arccos t,\arccos s)
$$
is the isotropic part of a positive definite kernel on $S^m\times S^M$.\ Further, if $g$ is nonconstant, then $f$ is the isotropic part of a strictly positive definite kernel on $S^m \times S^M$.
\end{prop}
\pf Consider the integral representation for $g$ as described above.\ If $x_1,\ldots,x_n\in S^m$, $w_1,\ldots,w_n\in S^M$ and $c_1,\ldots,c_n$ are real scalars, then
$$
\sum_{\mu,\nu=1}^n c_\mu c_\nu  f(x_\mu \cdot x_\nu, w_\mu \cdot w_\nu)  = \int_{[0,\infty)^2} \sum_{\mu,\nu=1}^n c_\mu c_\nu e^{-t\arccos(x_\mu\cdot x_\nu)-s\arccos(w_\mu\cdot w_\nu)}d\rho(t,s),
$$
that is,
$$
\sum_{\mu,\nu=1}^n c_\mu c_\nu  f(x_\mu \cdot x_\nu, w_\mu \cdot w_\nu)  = \int_{[0,\infty)^2} \sum_{\mu,\nu=1}^n c_\mu c_\nu e^{-t d_m(x_\mu\cdot x_\nu)-sd_M(w_\mu\cdot w_\nu)}d\rho(t,s),
$$
in which $d_m$ and $d_M$ are the usual geodesic distances on $S^m$ and $S^M$, respectively.\ A result proved in \cite{alexander} reveals that $d_m$ and $d_M$ are kernels of negative type.\ Consequently, the matrices with entries $-t d_m(x_\mu,x_\nu)-sd_M(w_\mu,w_\nu)$ is almost nonnegative definite (\cite[p.135]{don}).\ A classical result from the theory of positive definite kernels (\cite[p.74]{berg}) now implies that $\exp(-td_m-sd_M)$ is a positive definite kernel on $S^m \times S^M$.\ Thus, the initial quadratic form is nonnegative and the first assertion of the proposition is proved.\ As for the second one, it suffices to observe that if the points $(x_\mu,w_\mu)$ are distinct then the matrix with entries $-t d_m(x_\mu,x_\nu)-sd_M(w_\mu,w_\nu)$ has no pair of identical rows when $t,s>0$.\ In that case, the kernel $(x,z,y,w) \in (S^m \times S^M)^2 \to \exp[-td_m(x,y)-sd_M(z,w)]$ is, in fact, strictly positive definite on $S^m \times S^M$.\ If $g$ is nonconstant, then the original quadratic form is always positive unless all the $c_\mu$ are zero.\eop

To close the paper, we go the other way around, seeking positive definiteness on a single sphere from positive definiteness on a product of spheres.\ Two results in that direction are as follows.

\begin{prop} \label{restric} If $f$ is the continuous and isotropic part of a (strictly) positive definite kernel on $S^{m}\times S^M$, then $t  \to f(t,1)$ and $s \to f(1,s)$ are the isotropic parts of (strictly) positive definite kernels on $S^m$ and $S^M$ respectively.
\end{prop}

\begin{prop} \label{restric1} If $f$ is the continuous and isotropic part of a $DC$-strictly positive definite kernel on $S^{m}\times S^M$, then $t  \to f(t,t)$ is the isotropic part of a strictly positive definite kernel
 on $S^{m\wedge M}$, in which $m\wedge M=\min\{m,M\}$.
\end{prop}

We intend to provide solutions for some of the problems mentioned above in a subsequent paper.\ For now, we conclude this one mentioning a few relevant references that deals with similar questions on a single sphere: \cite{chen, cheney,gneiting,mene,ron,sun}.

%
%

\vspace*{2cm}

\noindent J. C. Guella, V. A. Menegatto, and A. P. Peron \\
Departamento de
Matem\'atica,\\ ICMC-USP - S\~ao Carlos, Caixa Postal 668,\\
13560-970 S\~ao Carlos SP, Brasil\\ e-mails: jeanguella@gmail.com; menegatt@icmc.usp.br; apperon@icmc.usp.br


\begin{thebibliography}{1}
\bibitem{alexander} Alexander, R.; Stolarsky, K. B., Extremal problems of distance geometry related to energy integrals.\ {\em Trans. Amer. Math. Soc.} 193 (1974), 1-31.
\bibitem{atkinson} Atkinson, K.; Han, Weimin, Spherical harmonics and approximations on the unit sphere: an introduction.\ Lecture Notes in Mathematics, 2044. Springer, Heidelberg, 2012.
\bibitem{berg}  Berg, C.; Christensen, J. P. R.; Ressel, P., Harmonic analysis on semigroups.\ Theory of positive definite and related functions.\ Graduate Texts in Mathematics, 100.\ Springer-Verlag, New York, 1984.
\bibitem{bochner} Bochner, S., Harmonic analysis and the theory of probability. University of California Press, Berkeley and Los Angeles, 1955.
\bibitem{chen} Chen, Debao; Menegatto, V. A.; Sun, Xingping, A necessary and sufficient condition for strictly positive definite functions on spheres.\ {\em Proc. Amer. Math. Soc.}
 131 (2003), no. 9, 2733-2740.
 \bibitem{cheney} Cheney, E. W., Approximation using positive definite functions.\ Approximation theory VIII, Vol. 1 (College Station, TX, 1995), 145-168, Ser. Approx. Decompos., 6, World Sci. Publ., River Edge, NJ, 1995.
\bibitem{dai} Dai, Feng; Xu, Yuan, Approximation theory and harmonic analysis on spheres and balls.\ Springer Monographs in Mathematics.\ Springer, New York, 2013.
\bibitem{don} Donoghue, W. F., Jr. Monotone matrix functions and analytic continuation. Die Grundlehren der mathematischen Wissenschaften, Band 207. Springer-Verlag, New York-Heidelberg, 1974.
\bibitem{limaye} Ghorpade, S. R.; Limaye, B. V., A course in multivariable calculus and analysis.\ Undergraduate Texts in Mathematics.\ Springer, New York, 2010.
\bibitem{gneiting} Gneiting, T., Strictly and non-strictly positive definite functions on spheres.\ {\em Bernoulli} 19 (2013), no. 4, 1327-1349.
\bibitem{groemer} Groemer, H., Geometric applications of Fourier series and spherical harmonics.\ Encyclopedia of Mathematics and its Applications, 61.\ Cambridge University Press, Cambridge, 1996.
\bibitem{horn} Horn, R. A.; Johnson, C. R., Matrix analysis.\ Second edition.\ Cambridge University Press, Cambridge, 2013.
\bibitem{mene0} Menegatto, V. A., Strict positive definiteness on spheres.\ {\em Analysis} 19 (1999), no. 3, 217-233.
\bibitem{mene} Menegatto, V. A.; Oliveira, C. P.; Peron, A. P., Strictly positive definite kernels on subsets of the complex plane.\ {\em Comput. Math. Appl.} 51 (2006), no. 8, 1233-1250.
\bibitem{samko} Miller, K. S.; Samko, S. G., Completely monotonic functions.\ {\em Integral Transform. Spec. Funct.} 12 (2001), no. 4, 389-402.
\bibitem{muller} M{\"u}ller, C., Analysis of spherical symmetries in Euclidean spaces.\ Applied Mathematical Sciences, 129.\ Springer-Verlag, New York, 1998.
\bibitem{ron} Ron, A.; Sun, Xingping, Strictly positive definite functions on spheres in Euclidean spaces.\ {\em Math. Comp.} 65 (1996), no. 216, 1513-1530.
\bibitem{schoen} Schoenberg, I. J.; Positive definite functions on spheres.\ {\em Duke Math. J.} 9, (1942), 96-108.
\bibitem{sun} Sun, Xingping, Strictly positive definite functions on the unit circle.\ {\em Math. Comp.} 74 (2005), no. 250, 709-721.
\bibitem{szego}
Szeg\"{o}, G., Orthogonal polynomials.\ Fourth edition.\ American Mathematical Society, Colloquium Publications, Vol. XXIII.\ American Mathematical Society, Providence, R.I., 1975.

\end{thebibliography}
\end{document}